\documentclass[11pt]{amsart}
\usepackage{amsmath}
\usepackage{amssymb}
\usepackage{amsfonts}

\newcommand{\C}{\mathbb{C}}
\newcommand{\Z}{\mathbb{Z}}

\newcommand{\cl}{{\rm cl}}
\newcommand{\intt}{{\rm int}}
\newcommand{\h}{{\mathcal H}}
\newcommand{\lh}{{\mathcal L}({\mathcal H})}
\newcommand{\ran}{{\rm ran}}

\DeclareMathSymbol{\subsetneqq}{\mathbin}{AMSb}{36}

\newtheorem{th1}{{\bf Theorem}}[section]
\newtheorem{thm}[th1]{{\bf Theorem}}
\newtheorem{lem}[th1]{{\bf Lemma}}
\newtheorem{prop}[th1]{{\bf Proposition}}
\newtheorem{cor}[th1]{{\bf Corollary}}

\theoremstyle{remark}
\newtheorem{rem}[th1]{Remark}

\theoremstyle{definition}

\newtheorem{exm}[th1]{Example}

\address{The Abdus Salam International Centre for Theoretical Physics,
Trieste, Italy}
\email{bourhim@ictp.trieste.it}

\begin{document}


\title{LOCAL SPECTRA OF OPERATOR WEIGHTED SHIFTS}
\author{A. Bourhim}
\date{}
\maketitle

\begin{abstract}
In this note, we study the local spectral properties of unilateral
operator weighted shifts.
\end{abstract}

\section{Introduction}
Let $\lh$ denote the algebra of all bounded linear operators
acting on a complex Hilbert space $\h$, and let ${\mathcal
A}:=(A_n)_{n\geq0}$ be a sequence of uniformly bounded invertible
operators of $\lh$. Let
$$\widehat{\h}=\sum\limits_{n=0}^{+\infty}\oplus\h_n,$$ where
$\h_n=\h$ for each $n\geq0$. It is a Hilbert space when equipped
with the inner product $$\langle
(x_n)_n,(y_n)_n\rangle_{\widehat{\h}}=\sum\limits_{n=0}^{+\infty}\langle
x_n,y_n\rangle_\h.$$ Therefore, the corresponding norm is given by
$$\|(x_n)_n\|_{\widehat{\h}}=\big(\sum\limits_n\|x_n\|_\h^2\big)^{\frac{1}{2}}.$$
The {\it unilateral operator weighted shift}, $S_u$, with the
weight sequence ${\mathcal A}=(A_n)_{n\geq0}$ is the operator on
$\widehat{\h}$ defined by
$$S_u(x_0,x_1,x_2,...)=(0,A_0x_0,A_1x_1,A_2x_2,...),~~((x_n)_n\in\widehat{\h}).$$

Operator weighted shifts were first introduced by A. Lambert
\cite{lambert}, and have been studied by many authors (see for
example \cite{ben}, \cite{herrero}, \cite{jue}, \cite{lan}, and
\cite{you}). In the case when $\dim\h=1$, they are exactly the
scalar weighted shifts which have been widely studied. An
excellent survey of the investigation of the spectral theory of
such operators was given by A. L. Shields \cite{shi}. Moreover,
several known results for the scalar case have been generalized
and extended to the setting of operator weighted shifts. However,
the question of determining the local spectral properties for
operator weighted shifts is natural and has been initiated in
\cite{zguitti}. While, the investigation of these properties for
scalar weighted shifts has been studied in \cite{bourhim} and
\cite{mmn}. The main goal of the present note is to
study and examine whether or not the results obtained in
\cite{bourhim} remain valid for unilateral operator weighted
shifts. We give necessary and sufficient conditions for a
unilateral operator weighted shift to satisfy Dunford's condition
$(C)$ or Bishop's property $(\beta)$. Unlike the scalar weighted
shift operators, we show that they are examples of unilateral
operator weighted shifts possessing Bishop's property $(\beta)$
with large approximate point spectrum and without fat local
spectra.

For an operator $T\in \lh$, let, as usual, $T^*$, $\sigma(T)$,
$\sigma_{ap}(T)$, $\sigma_p(T)$, and $r(T)$ denote the adjoint,
the spectrum, the approximate point spectrum, the point spectrum,
and the spectral radius of $T$, respectively. Let
$m(T):=\inf\{\|Tx\|:\|x\|=1\}$ denote the lower bound of $T$, and
note that the sequence $\big(m(T^n)^\frac{1}{n}\big)_{n\geq1}$
converges and its limit, denoted $r_1(T)$, equals
$\sup\limits_{n\geq 1}[m(T^n)]^\frac{1}{n}$ (see \cite{zemanek}).
Let $T\in\lh$; for an element $x\in\h$, let $\sigma_{_T}(x)$,
$\rho_{_T}(x):=\C\backslash\sigma_{_T}(x)$, and
$r_T(x):=\limsup\limits_{n\to+\infty}\|T^nx\|^{\frac{1}{n}}$ be
the local spectrum, the local resolvent set and the local spectral
radius of $T$ at $x$ , respectively (see \cite{cf} and \cite{ln}).
The operator $T$ is said to have the {\it single--valued extension
property} at a complex number $\lambda_0\in\C$ if for every open
disc $U$ centered at $\lambda_0$, the only analytic solution of
the equation $(T-\lambda)f(\lambda )=0,~(\lambda\in U)$ is the
zero function $f\equiv 0$. Denote by $\Re(T)$ the set of all
complex numbers on which $T$ fails to have the single--valued
extension property and recall that $T$ is said to have the {\it
single--valued extension property} provided that $\Re(T)$ is
empty. The reader is reminded that in the case $T$ has the
single--valued extension property, the {\it local resolvent} of
$x$ is the unique analytic $\h-$valued function,
$\widetilde{x}(.)$, satisfying
$(T-\lambda)\widetilde{x}(\lambda)=x,~(\lambda\in\rho_{_T}(x))$.
Also, recall that an operator $T\in\h$ is said to satisfy {\it
Dunford's condition $(C)$} provided that for every closed subset
$F$ of $\C$, the linear subspace,
$$\h_{_T}(F):=\{x\in\h:\sigma_{_T}(x)\subset F\}$$ is closed.
Moreover, $T$ is said to have {\it fat local spectra} if
$\sigma_{_T}(x)=\sigma(T)$ for all non-zero $x\in\h$. It is well
known that every operator which satisfies Dunford's condition
$(C)$ has the single--valued extension property and it turns out
that Dunford's condition $(C)$ follows from fat local spectra
property.

Throughout this note, let $S_u$ be a unilateral operator weighted
shift with weight sequence ${\mathcal A}:=(A_n)_{n\geq0}$, and let
$(B_n)_{n\geq0}$ be the sequence given by $$ B_n=\left\{
\begin{array}{lll}
A_{n-1}A_{n-2}...A_1A_0&\mbox{if }n>0\\
\\
1&\mbox{if }n=0
\end{array}
\right. $$ Define
$$r_2(S_u):=\frac{1}{\limsup\limits_{n\to+\infty}\|{B_n}^{-1}\|^{\frac{1}{n}}},~~
r_3(S_u):=\frac{1}{\liminf\limits_{n\to+\infty}\|{B_n}^{-1}\|^{\frac{1}{n}}},$$
$$R_2^+(S_u):=\sup\limits_{x\in\h,~x\not=0}
\bigg\{\frac{1}{\limsup\limits_{n\to+\infty}\|{B_n^*}^{-1}x\|^{\frac{1}{n}}}\bigg\}=
\sup\limits_{x\in\h,~\|x\|=1}
\bigg\{\frac{1}{\limsup\limits_{n\to+\infty}\|{B_n^*}^{-1}x\|^{\frac{1}{n}}}\bigg\},$$
$$R_2^-(S_u):=\inf\limits_{x\in\h,~x\not=0}
\bigg\{\frac{1}{\limsup\limits_{n\to+\infty}\|{B_n^*}^{-1}x\|^{\frac{1}{n}}}\bigg\}=
\inf\limits_{x\in\h,~\|x\|=1}
\bigg\{\frac{1}{\limsup\limits_{n\to+\infty}\|{B_n^*}^{-1}x\|^{\frac{1}{n}}}\bigg\},$$
$$R_3^+(S_u):=\sup\limits_{x\in\h,~x\not=0}\bigg\{\limsup\limits_{n\to+\infty}
\|B_nx\|^{\frac{1}{n}}\bigg\}=\sup\limits_{x\in\h,~\|x\|=1}\bigg\{\limsup\limits_{n\to+\infty}
\|B_nx\|^{\frac{1}{n}}\bigg\},$$ and
$$R_3^-(S_u):=\inf\limits_{x\in\h,~x\not=0}\bigg\{\limsup\limits_{n\to+\infty}
\|B_nx\|^{\frac{1}{n}}\bigg\}=\inf\limits_{x\in\h,~\|x\|=1}\bigg\{\limsup\limits_{n\to+\infty}
\|B_nx\|^{\frac{1}{n}}\bigg\}.$$ Note that $$r_1(S_u)\leq
r_2(S_u)\leq R_2^-(S_u)\leq R_2^+(S_u),$$ and $$r_3(S_u)\leq
R_3^-(S_u)\leq R_3^+(S_u)\leq r(S_u).$$ Note also that for a
scalar weighted shift $S_u$, we have $$r_1(S_u)\leq r_2(S_u)=
R_2^-(S_u)= R_2^+(S_u)\leq r_3(S_u)= R_3^-(S_u)=R_3^+(S_u)\leq
r(S_u).$$

Finally, we would like to record and without further mention a
notation that we will use repeatedly throughout this note. For
every $x\in\h$, we write $$x^{(n)}=(0,...,0,x,0,...),~~(n\geq0)$$
for the element of ${\widehat{\h}}$ for which all the coordinates
are zero except the $n$th coordinate which is equal $x$, and note
that
\begin{equation}\label{localradius}
r_{S_u}(x^{(k)})=\limsup\limits_{n\to+\infty}\|B_{n+k}B_k^{-1}x\|^{\frac{1}{n}}.
\end{equation}

\section{Preliminaries and elementary background}
In this section, we assemble some elementary results that are very
much on the straightforward side and therefore the proofs will be
omitted.

\begin{prop}\label{straightforward1}
Assume that $T\in\lh$ is an operator for which
$\bigcap\limits_{n\geq0}T^n\h=\{0\}$. The following statements
hold.
\begin{itemize}
\item[$(a)$]$\{\lambda\in\C:|\lambda|\leq
r_1(T)\}\subset\sigma_{_T}(x)$ for every non-zero element
$x\in\h$.
\item[$(b)$]$\sigma_p(T)\subset\{0\}$.
\item[$(c)$]Each $\sigma_{_T}(x)$ is connected.
\item[$(d)$]$\sigma(T)$ is a connected set and satisfies
$\{\lambda\in\C:|\lambda|\leq r_1(T)\}\subset\sigma(T).$ In
particular, if $\sigma(T)$ is circularly symmetric about the
origin, then $$\sigma(T)=\{\lambda\in\C:|\lambda|\leq r(T)\}.$$
\end{itemize}
\end{prop}

Evidently, the unilateral operator weighted shift $S_u$ satisfies
the condition that
$\bigcap\limits_{n\geq0}S_u^n{\widehat{\h}}=\{0\}$, and its
spectrum is rotationally symmetric. Therefore, the next result is
an immediate consequence of proposition \ref{straightforward1}.

\begin{cor}\label{straightforward2}The following statements hold.
\begin{itemize}
\item[$(a)$]For every non-zero element
$x\in{\widehat{\h}}$, the local spectrum, $\sigma_{_{S_u}}(x)$, of
$S_u$ at $x$ is connected and satisfies
$\{\lambda\in\C:|\lambda|\leq
r_1(S_u)\}\subset\sigma_{_{S_u}}(x).$
\item[$(b)$]The spectrum of $S_u$ is the disc $\{\lambda\in\C:|\lambda|\leq r(S_u)\}.$
\end{itemize}
\end{cor}

\begin{prop}For every $n\geq1$, we have
$$\|S_u^n\|=\sup\limits_{k\geq0}\|B_{n+k}B_k^{-1}\|,\mbox{ and }
m(S_u^n)=\inf\limits_{k\geq0}\big\{\frac{1}{\|B_{k}B_{n+k}^{-1}\|}\big\}.$$
Thus, $$r(S_u)=\lim\limits_{n\to+\infty}\bigg[\sup\limits_{k\geq0}
\|B_{n+k}B_k^{-1}\|\bigg]^{\frac{1}{n}},\mbox{ and
}r_1(S_u)=\lim\limits_{n\to+\infty}\bigg[
\inf\limits_{k\geq0}\big\{\frac{1}{\|B_{k}B_{n+k}^{-1}\|}\big\}\bigg]^{\frac{1}{n}}.$$
\end{prop}

\begin{prop}The adjoint of $S_u$ is given by
$$S_u^*x=
(A_0^*x_1,A_1^*x_2,A_2^*x_3,...),~(x=(x_0,x_1,...)\in{\widehat{\h}}).$$
\end{prop}

\section{Local spectra of $S_u$}

We begin this section with the following result that gives a
necessary and sufficient condition for $S_u^*$ to enjoy the
single--valued extension property.

\begin{lem}\label{spectreponctuel}The following statments hold.
\begin{itemize}
\item[$(a)$] $\sigma_p({S_u})=\emptyset$.
\item[$(b)$] $\{0\}\cup\{\lambda\in\C:|\lambda|<R_2^+(S_u)\}
\subset\sigma_p({S_u^*})\subset\{\lambda\in\C:|\lambda|\leq
R_2^+(S_u)\}.$
\item[$(c)$] $S_u^*$ has the single--valued
extension property if and only if $R_2^+(S_u)=0$. Moreover, we
always have $$\Re(S_u^*)=\{\lambda\in\C:|\lambda|<R_2^+(S_u)\}.$$
\end{itemize}
\end{lem}
\begin{proof}
$(a)$ By proposition \ref{straightforward1}-$(b)$, we have
$\sigma_p({S_u})\subset\{0\}$. As $S_u$ is injective, we note that
$\sigma_p({S_u})=\emptyset$.

$(b)$ Suppose that $\lambda\in\C$ is an eigenvalue for $S_u^*$ and
that $(x_n)_n$ is a corresponding eigenvector. We have
$$(A_0^*x_1,A_1^*x_2,A_2^*x_3,...)=(\lambda x_0,\lambda
x_1,\lambda x_2,...).$$ This shows that
$$x_n=\lambda^n{B_n^*}^{-1}x_0,~(n\geq0).$$ Therefore,
$$\|x\|^2=\sum\limits_{n\geq0}|\lambda|^{2n}\|{B_n^*}^{-1}x_0\|^2.$$
By the Cauchy-Hadamard formula for the radius of convergence, we
get that $$|\lambda|\leq
\frac{1}{\limsup\limits_{n\to+\infty}\|{B_n^*}^{-1}x_0\|^{\frac{1}{n}}}.$$
Thus, $$\sigma_p({S_u^*})\subset\{\lambda\in\C:|\lambda|\leq
R_2^+(S_u)\}.$$

Now, let us prove that
$$\{0\}\cup\{\lambda\in\C:|\lambda|<R_2^+(S_u)\}
\subset\sigma_p({S_u^*}).$$ It is clear that for every $x\in\h$,
we have $S_u^*x^{(0)}=0$; hence, $0\in\sigma_p({S_u^*})$. If
$R_2^+(S_u)=0$, then there is nothing to prove; thus, we may
assume that $R_2^+(S_u)>0$. Let $\lambda\in\C$ such that
$|\lambda|<R_2^+(S_u)$. So, there is a non-zero $x_0\in\h$ such
that
$|\lambda|<\frac{1}{\limsup\limits_{n\to+\infty}\|{B_n^*}^{-1}x_0\|^{\frac{1}{n}}}.$
We have $\big(S_u^*-\lambda\big)k_{x_0}(\lambda)=0,$ where
$k_{x_0}(\lambda)=\sum\limits_{n\geq0}\oplus
\lambda^n{B_n^*}^{-1}x_0.$ This shows that
$$\{\lambda\in\C:|\lambda|<R_2^+(S_u)\}
\subset\sigma_p({S_u^*}),$$ and the desired statement holds.

 $(c)$ In view of the statement $(b)$ and the fact that
$\Re(S_u^*)\subset\intt\big(\sigma_p({S_u^*})\big)$, we have
$\Re(S_u^*)\subset\{\lambda\in\C:|\lambda|<R_2^+(S_u)\}.$

Conversely, let $x$ be a non-zero element of $\h$ and set
$$U_x:=\big\{\lambda\in\C:|\lambda|<\frac{1}{{\limsup\limits_{n\to+\infty}
\|{B_n^*}^{-1}x\|^{\frac{1}{n}}}}\big\},$$and $$
k_{x}(\lambda):=\sum\limits_{n\geq0}\oplus
\lambda^n{B_n^*}^{-1}x,~(\lambda\in U_x).$$ Since
$\big(S_u^*-\lambda\big)k_{x}(\lambda)=0,~\mbox{ for all
}\lambda\in U_x,$ and $x$ is an arbitrary non-zero element of
$\h$, we have
$$\{\lambda\in\C:|\lambda|<R_2^+(S_u)\}=\bigcup\limits_{x\in\h,~x\not=0}U_x\subset\Re(S_u^*).$$
The proof is therefore complete.
\end{proof}

The following result refine the local spectral inclusion given in
corollary \ref{straightforward2}.
\begin{prop}\label{fatlocal}For every non-zero $y=(y_0,y_1,y_2,...)\in{\widehat{\h}}$, we have
$$\{\lambda\in\C:|\lambda|\leq
R_2^-(S_u)\}\subset\sigma_{_{S_u}}(y).$$ In particular, if
$r(S_u)=R_2^-(S_u)$ then $S_u$ has fat local spectra.
\end{prop}
\begin{proof}
As $\bigcap\limits_{n\geq 0}S_u^n{\widehat{\h}}=\{0\}$, we have
$0\in\sigma_{_{S_u}}(y).$ Thus, we may assume that $R_2^-(S_u)>0$.
Let $O:=\{\lambda\in\C:|\lambda|<R_2^-(S_u)\}$, and let $x$ be a
non-zero element of $\h$. Consider the following analytic
${\widehat{\h}}$--valued function on $O$,
$$k_x(\lambda)=\sum\limits_{n\geq0}\oplus
\lambda^n{B_n^*}^{-1}x.$$ We have
$(S_u-\lambda)^*k_x(\overline{\lambda})=0$ for every $\lambda\in
O$. Now, let $y=(y_0,y_1,y_2,...)\in{\widehat{\h}}$ such that
$O\cap\rho_{_{S_u}}(y)\not=\emptyset$. So, for every $\lambda\in
O\cap\rho_{_{S_u}}(y),$ we have
\begin{eqnarray*}
\sum\limits_{n\geq0}\langle
y_n,{B_n^*}^{-1}x\rangle_\h\lambda^n&=& \langle
y,k_x(\overline{\lambda})\rangle_{\widehat{\h}}\\
&=&\langle(S_u-\lambda)\widetilde{y}(\lambda),k_x(\overline{\lambda})\rangle_{\widehat{\h}}\\
&=&\langle\widetilde{y}(\lambda),(S_u-\lambda)^*k_x(\overline{\lambda})\rangle_{\widehat{\h}}\\
&=&0.
\end{eqnarray*}
Hence, for every $n\geq0$, we have $$\langle
y_n,{B_n^*}^{-1}x\rangle_\h=0.$$ Since $x$ is an arbitrary element
of $\h$, we have $y=0$; and the proof is complete.
\end{proof}

In view of proposition \ref{fatlocal}, we note that
$R_2^-(S_u)\leq r_{S_u}(x)$, for all non-zero
$x=(x_0,x_1,x_2,...)\in{\widehat{\h}}$. The following gives more
information about local spectral radii of $S_u$.

\begin{prop}\label{lradius}
For every non-zero element $x=(x_0,x_1,...)\in{\widehat{\h}}$, we
have $$R_3^-(S_u)\leq r_{S_u}(x)\leq r(S_u).$$ Moreover, if
$x=(x_0,x_1,...)$ is a non-zero finitely supported element of
${\widehat{\h}}$, then
\begin{equation}\label{dense}
R_3^-(S_u)\leq
r_{S_u}(x)=\max\limits_{k\geq0}\big(r_{S_u}(x_k^{(k)})\big)\leq
R_3^+(S_u).
\end{equation}
\end{prop}
\begin{proof}Let $x=(x_0,x_1,...)$ be a non-zero element of
${\widehat{\h}}$; so, there is an integer $k_0\geq0$ such that
$x_{k_0}\not=0$. Since,
$$\|S_u^nx\|^2=\sum\limits_{k=0}^{+\infty}\|B_{n+k}B_k^{-1}x_k\|^2,~\forall
n\geq0,$$ we have
$$\|B_{n+k_0}B_{k_0}^{-1}x_{k_0}\|^{\frac{1}{n+k_0}}\leq\|S_u^nx\|^{\frac{1}{n+k_0}},~\forall
n\geq0 .$$ Now, taking $\limsup$ as $n\to+\infty$, we get
$$R_3^-(S_u)\leq
\limsup\limits_{n\to+\infty}\|B_{n+k_0}B_{k_0}^{-1}x_{k_0}\|^{\frac{1}{n+k_0}}\leq
r_{S_u}(x),$$ as desired.

$(b)$ Assume that $x=(x_0,x_1,...)$ is a non-zero finitely
supported element of ${\widehat{\h}}$. As above, we have
$$\|B_{n+k}B_k^{-1}x_k\|^{\frac{1}{n}}\leq\|S_u^nx\|^{\frac{1}{n}},~\forall
n,~k\geq0 .$$ By taking $\limsup$ as $n\to+\infty$, we get
$r_{S_u}(x_k^{(k)})\leq r_{S_u}(x),~\forall k\geq0$. Hence,
$$\max\limits_{k\geq0}\big(r_{S_u}(x_k^{(k)})\big)\leq
r_{S_u}(x).$$ As
$\sigma_{_{S_u}}(x)\subset\bigcup\limits_{k\geq0}\sigma_{_{S_u}}(x_k^{(k)})$,
and $r_{S_u}(y)=\max\{|\lambda|:\lambda\in\sigma_{_{S_u}}(y)\}$
for every non-zero $y\in{\widehat{\h}}$, we obtain
$r_{S_u}(x)\leq\max\limits_{k\geq0}\big(r_{S_u}(x_k^{(k)})\big).$
Hence,
$$r_{S_u}(x)=\max\limits_{k\geq0}\big(r_{S_u}(x_k^{(k)})\big).$$
On the other hand, we have
$r_{S_u}(x_k^{(k)})=r_{S_u}((B_k^{-1}x_k)^{(0)}),~\forall k\geq0$.
This shows that
$$r_{S_u}(x)=\max\limits_{k\geq0}\big(r_{S_u}(x_k^{(k)})\big)\leq
R_3^+(S_u).$$ Therefore, the desired result holds.
\end{proof}

For every $x=(x_0,x_1,...)\in{\widehat{\h}}$, we set
$$R_{{\mathcal
A}}(x):=\frac{1}{\limsup\limits_{n\to+\infty}\|B_n^{-1}x_n\|^{\frac{1}{n}}}.$$
Obviously, if $x$ is a non-zero element of $\widehat{\h}$, then
$r_2(S_u)\leq R_{{\mathcal A}}(x)\leq+\infty$.

\begin{thm}\label{localspectrum}
For every non-zero element $x=(x_0,x_1,...)\in{\widehat{\h}}$, we
have $$\{\lambda\in\C:|\lambda|\leq \min\big(R_{{\mathcal
A}}(x),r_3(S_u)\big) \}\subset\sigma_{_{S_u}}(x).$$ Moreover, if
$x=(x_0,x_1,...)$ is a non-zero finitely supported element of
${\widehat{\h}}$, then $$\{\lambda\in\C:|\lambda|\leq R_3^-(S_u)
\}\subset\sigma_{_{S_u}}(x).$$
\end{thm}
\begin{proof}
Let $x=(x_0,x_1,...)$ be a non-zero element of ${\widehat{\h}}$.
If $\min\big(R_{{\mathcal A}}(x),r_3(S_u)\big)=0$, then there is
nothing to prove since $0\in\sigma_{_{S_u}}(x)$. Thus we may
suppose that $\min\big(R_{{\mathcal A}}(x),r_3(S_u)\big)>0$. Now,
for each $n\geq0$, let
$$F_n(\lambda)=-\frac{B_nx_0}{\lambda^{n+1}}-\frac{B_nB_1^{-1}x_1}{\lambda^{n}}-
\frac{B_nB_2^{-1}x_2}{\lambda^{n-1}}-...-\frac{x_n}{\lambda}
,~(\lambda\in\C\backslash\{0\}),$$ and $$G_n(\lambda)=x_0+\lambda
B_1^{-1}x_1+\lambda^2
B_2^{-1}x_2+...+\lambda^nB_n^{-1}x_n,~(\lambda\in\C).$$ We have,
\begin{equation}\label{AB}
F_n(\lambda)=\frac{-1}{\lambda^{n+1}}B_nG_n(\lambda),~(\lambda\in\C\backslash\{0\}).
\end{equation}
By writing $\widetilde{x}(\lambda):=(f_0(\lambda),
f_1(\lambda),f_2(\lambda),...),~\lambda\in\rho_{_{S_u}}(x)$, we
get from the equation,
$$(S_u-\lambda)\widetilde{x}(\lambda)=x,~\lambda\in\rho_{_{S_u}}(x),$$
that for every $\lambda\in\rho_{_{S_u}}(x)$, we have $$ \left\{
\begin{array}{lll}
-\lambda f_0(\lambda)=x_0\\
\\
A_nf_n(\lambda)-\lambda f_{n+1}(\lambda)=x_{n+1}&\mbox{ for every
}n\geq 0.
\end{array}
\right. $$ Therefore, for every $n\geq0$ and for every
$\lambda\in\rho_{_{S_u}}(x)$, we have
\begin{eqnarray*}
f_n(\lambda)&=&-\frac{B_nx_0}{\lambda^{n+1}}-\frac{B_nB_1^{-1}x_1}{\lambda^{n}}-
\frac{B_nB_2^{-1}x_2}{\lambda^{n-1}}-...-\frac{x_n}{\lambda}\\
&=&F_n(\lambda).
\end{eqnarray*}
Since
$\|\widetilde{x}(\lambda)\|^2=\sum\limits_{n\geq0}\|f_n(\lambda)\|^2<+\infty$
for every $\lambda\in\rho_{_{S_u}}(x)$, it then follows that
\begin{equation}\label{lim}
\lim\limits_{n\to+\infty}F_n(\lambda)=\lim\limits_{n\to+\infty}f_n(\lambda)=0\mbox{
for every }\lambda\in\rho_{_{S_u}}(x).
\end{equation}
We shall show that (\ref{lim}) is not satisfied for most of the
points in the open disc
$V(x):=\{\lambda\in\C:|\lambda|<\min\big(R_{{\mathcal
A}}(x),r_3(S_u)\big)\}$. It is clear that the sequence
$(G_n)_{n\geq0}$ converges uniformly on compact subsets of $V(x)$
to the non-zero power series
$G(\lambda)=\sum\limits_{n\geq0}\lambda^nB_n^{-1}x_n$. Now, let
$\lambda_0\in V(x)\backslash\{0\}$ such that $G(\lambda_0)\not=0$;
there is $\epsilon>0$ and an integer $n_0$ such that
$\epsilon<\|G_n(\lambda_0)\|$ for every $n\geq n_0$. On the other
hand, $|\lambda_0|<r_3(S_u)$, then there is a subsequence
$(n_k)_{k\geq0}$ of integers greater than $n_0$ such that
$|\lambda_0|^{n_k}\|B_{n_k}^{-1}\|<1$. Thus, it follows from
(\ref{AB}) that for every $k\geq0$, we have
\begin{eqnarray*}
\|F_{n_k}(\lambda_0)\|&=&|\frac{-1}{\lambda_0^{n_k+1}}|\|B_{n_k}G_{n_k}(\lambda_0)\|\\
&\geq&\frac{1}{|\lambda_0^{n_k+1}|\|B_{n_k}^{-1}\|}\|G_{n_k}(\lambda_0)\|\\
&\geq&\frac{\epsilon}{|\lambda_0|}.
\end{eqnarray*}
And so, by (\ref{lim}), $\lambda_0\not\in\rho_{_{S_u}}(x)$. Since
the set of zeros of $G$ is at most countable, we have
$\{\lambda\in\C:|\lambda|\leq\min\big(R_{{\mathcal
A}}(x),r_3(S_u)\big)\}\subset\sigma_{_{S_u}}(x).$

Now, assume that $x=(x_0,x_1,...)$ is a non-zero finitely
supported element of ${\widehat{\h}}$, and $k_0$ is the largest
integer $n\geq0$ for which $x_n\not=0$. Conserve the same
notations as above and note that, for every $n\geq k_0$, we have
$$F_n(\lambda)=\frac{-1}{\lambda^{n+1}}B_nG(\lambda),~(\lambda\in\C\backslash\{0\}),$$
where $$G(\lambda):=x_0+\lambda B_1^{-1}x_1+\lambda^2
B_2^{-1}x_2+...+\lambda^{k_0}B_{k_0}^{-1}x_{k_0},~(\lambda\in\C).$$
Let $W(x):=\{\lambda\in\C:|\lambda|<R_3^-(S_u)\}$, and let
$\lambda_0\in W(x)\backslash\{0\}$ such that $G(\lambda_0)\not=0$.
As $|\lambda_0|<R_3^-(S_u)\leq
\limsup\limits_{n\to+\infty}\|B_nG(\lambda_0)\|^{\frac{1}{n}}$, we
note that the series $\sum\limits_{n\geq0}\|F_n(\lambda_0)\|^2$
diverges. Hence, $\lambda_0\in\sigma_{_{S_u}}(x)$, and therefore,
$$\{\lambda\in\C:|\lambda|\leq
R_3^-(S_u)\}\subset\sigma_{_{S_u}}(x).$$
\end{proof}

For every $x\in\h$,
we write
$${\widehat{\h}}(x):=\bigvee\{\big(B_nx\big)^{(n)}:n\geq0\},$$
where "$\bigvee$" denotes the closed linear span. It is shown in proposition
4.3.5 of \cite{zguitti} that for every non-zero $x\in\h$, we have
$$\sigma_{_{S_u}}(x^{(n)}) =\{\lambda\in\C:|\lambda|\leq
r_{S_u}(x^{(n)})\},~~(n\geq0).$$
We refine this result as follows; our proof is inspired by an argument of \cite{bourhim}.

\begin{prop}Let $x$ be a non-zero element of $\h$, and let
$y\in{\widehat{\h}}(x)$. The following statements hold.
\begin{itemize}
\item[$(a)$]If $R_{{\mathcal A}}(y)>r_{S_u}(x^{(0)})$, then
$\sigma_{_{S_u}}(y)=\{\lambda\in\C:|\lambda|\leq
r_{S_u}(x^{(0)})\}.$
\item[$(b)$]If $R_{{\mathcal A}}(y)\leq r_{S_u}(x^{(0)})$, then $\{\lambda\in\C:|\lambda|\leq
R_{{\mathcal A}}(y)\}\subset\sigma_{_{S_u}}(y).$
\end{itemize}
\end{prop}
\begin{proof}Let $x$ be a non-zero element of $\h$, and let us
first show that $$\sigma_{_{S_u}}(x^{(0)})
=\{\lambda\in\C:|\lambda|\leq r_{S_u}(x^{(0)})\}.$$ To do this it
suffices to prove that $\{\lambda\in\C:|\lambda|\leq
r_{S_u}(x^{(0)})\}\subset\sigma_{_{S_u}}(x^{(0)}).$ Indeed, as in
the proof of theorem \ref{localspectrum}, we trivially have
$$\widetilde{x^{(0)}}(\lambda)=(-\frac{x}{\lambda},-\frac{B_1x}{\lambda^2},
-\frac{B_2x}{\lambda^3},...),~(
\lambda\in\rho_{_{S_u}}(x^{(0)})).$$ In particular, we have
$\|\widetilde{x^{(0)}}(\lambda)\|_{\widehat{\h}}^2=
\sum\limits_{k=0}^{+\infty}\frac{\|B_kx\|_\h^2}{|\lambda|^{2(k+1)}},~(
\lambda\in\rho_{_{S_u}}(x^{(0)})).$ This implies that
$\rho_{_{S_u}}(x^{(0)})\subset\{\lambda\in\C:r_{S_u}(x^{(0)})\leq
|\lambda|\}.$ Or, equivalently, $$\{\lambda\in\C:|\lambda|<
r_{S_u}(x^{(0)})\}\subset\sigma_{_{S_u}}(x^{(0)}).$$ As
$\sigma_{_{S_u}}(x^{(0)})$ is a closed set, the desired identity
holds.

$(a)$ Assume that $y=\sum\limits_{n=0}^{+\infty}a_n\big(B_nx\big)^{(n)}$
is a non-zero element of ${\widehat{\h}}(x)$ for which
$R_{{\mathcal A}}(y)>r_{S_u}(x^{(0)})$. In this case the function
$f(\lambda):=\sum\limits_{n\geq
0}a_n\lambda^n$ is analytic on the open disc
$\{\lambda\in\C:|\lambda|<R_{{\mathcal A}}(y)\}$ which is a neighborhood of
$\sigma_{_{S_u}}(x^{(0)})$. Let $r$ be a real number such that
$r_{S_u}(x^{(0)})<r<R_{{\mathcal A}}(y)$, we have
\begin{eqnarray*}
f(S_u,x^{(0)})&:=&\frac{-1}{2\pi
i}\oint_{|\lambda|=r}f(\lambda)\widetilde{x^{(0)}}(\lambda)d\lambda\\
&=&\frac{-1}{2\pi
i}\oint_{|\lambda|=r}f(\lambda)\bigg(-\sum\limits_{n\geq
0}\frac{S_u^nx^{(0)}}{\lambda^{n+1}}\bigg)d\lambda\\
&=&y.
\end{eqnarray*}
And so, by theorem 2.12 of \cite{williams}, we have
$$\sigma_{_{S_u}}(y)=\sigma_{_{S_u}}(f(S_u,x^{(0)}))=\sigma_{_{S_u}}(x^{(0)})=
\{\lambda\in\C:|\lambda|\leq r_{S_u}(x^{(0)})\}.$$

$(b)$ The proof of the second statement is similar to the one of theorem
\ref{localspectrum} if, for every integer $n\geq0$, we take
$$F_n(\lambda):=-\big(\frac{a_n}{\lambda^{n+1}}+\frac{a_1}{\lambda^{n}}+
\frac{a_2}{\lambda^{n-1}}+...+\frac{a_n}{\lambda}\big)B_nx
,~(\lambda\in\C\backslash\{0\}),$$ and
$$G_n(\lambda):=a_0+a_1\lambda
+a_2\lambda^2+...+a_n\lambda^n,~(\lambda\in\C).$$
\end{proof}

\section{Dunford's condition $(C)$ and Bishop's property $(\beta)$ for $S_u$}

Before outlining the statement of the main results of this
section, let us recall a few more notions and properties from the
local spectral theory which will be needed in the sequel. An
operator $T\in\lh$ is said to be {\it hyponormal} if
$\|T^*x\|\leq\|Tx\|$ for all $x\in\h$. It is said be {\it
subnormal} if it has a normal extension, this means that there is
a normal operator $N$ on a Hilbert space ${\mathcal K}$,
containing $\h$, such that $\h$ is a closed invariant subspace of
$N$ and the restriction $N_{|\h}$ coincides with $T$. Note that
every subnormal operator is hyponormal, but the converse is false
(see \cite{con}). For an open subset $U$ of $\C$, let ${\mathcal
O}(U,\h)$ denote as usual the Fr\'echet space of all analytic
$\h-$valued functions on $U$. An operator $T\in\lh$ is said to
possess Bishop's property $(\beta)$ if the continuous mapping
\begin{eqnarray*}
T_U&:&{\mathcal O}(U,\h)\longrightarrow{\mathcal O}(U,\h)\\
&&\hspace{1.15cm}f\longmapsto (T-z)f
\end{eqnarray*}
is injective with closed range for each open subset $U$ of $\C$.
It is known that hyponormal operators possess Bishop's property
$(\beta)$ (see \cite{putinar}) and it turns out that Dunford's
condition $(C)$ follows from Bishop's property $(\beta)$. Let
$\lambda_0\in\C$; recall that an operator $T\in \lh$ is said to
possess {\it Bishop's property $(\beta)$} at $\lambda_0$ if there
is an open neighbourhood $V$ of $\lambda_0$ such that for every
open subset $U$ of $V$, the mapping $T_{_U}$ is injective and has
a closed range. Note that if $T$ possesses Bishop's property
$(\beta)$ at any point $\lambda\in\C$ then $T$ possesses Bishop's
classical property $(\beta)$. Finally, for any operator $T\in\lh$,
we shall denote $$\sigma_{\beta}(T):=\big\{\lambda\in\C:T\mbox{
fails to possess Bishop's property }(\beta)\mbox{ at
}\lambda\big\}.$$ It is a closed subset of $\sigma_{ap}(T)$ (see
proposition 2.1 of \cite{bourhim}).

The following result gives necessary conditions for the operator
weighted shift, $S_u$, to enjoy Dunford's condition $(C)$.

\begin{thm}\label{dcc}
If $S_u$ satisfies Dunford's condition $(C)$, then
$r(S_u)=R_3^+(S_u)$. Moreover, for every non-zero $x\in\h$, we
have
\begin{equation}\label{ddcc}
\limsup\limits_{n\to+\infty}\|B_nx\|^{\frac{1}{n}}=
\lim\limits_{n\to+\infty}
\bigg[\sup\limits_{k\geq0}\frac{\|B_{n+k}x\|}{\|B_kx\|}\bigg]^{\frac{1}{n}}.
\end{equation}
\end{thm}
\begin{proof}To prove $R_3^+(S_u)=r(S_u)$, it suffices to show that
$r(S_u)\leq R_3^+(S_u)$. Since each $B_k$ is an invertible
operator, we note that
$$R_3^+(S_u)=\sup\limits_{x\in\h,~x\not=0}\big(r_{S_u}(x^{(k)})\big),~\forall
k\geq0.$$ Now, assume that $S_u$ satisfies Dunford's condition
$(C)$, and let $$F:=\{\lambda\in\C:|\lambda|\leq R_3^+(S_u)\}.$$
It follows from (\ref{dense}) that ${\widehat{\h}}_{_{S_u}}(F)$
contains a dense subspace of ${\widehat{\h}}$. As the subspace
${\widehat{\h}}_{_{S_u}}(F)$ is closed, we have
${\widehat{\h}}_{_{S_u}}(F)=\h$; therefore,
$\sigma_{_{S_u}}(x)\subset F$ for every $x\in{\widehat{\h}}$. And
so,
$\sigma({S_u})=\bigcup\limits_{x\in\h}\sigma_{_{S_u}}(x)\subset F$
(see proposition 1.3.2 of \cite{ln}). Hence, $r(S_u)\leq
R_3^+(S_u)$, as desired.

Let $x$ be a non-zero element of $\h$ and let us now establish the
identity (\ref{ddcc}). Since $S_u$ satisfies Dunford's condition
$(C)$, we note that $S_u$ restricted to ${\widehat{\h}}(x)$
satisfies also Dunford's condition $(C)$ (see proposition 1.2.21
of \cite{ln}). Now, note that $(v_n)_{n\geq0}$ is an orthonormal
basis of ${\widehat{\h}}(x)$, where
$$v_n:=\frac{(B_nx)^{(n)}}{\|B_nx\|},~(n\geq0).$$ We have
$$S_uv_n=\frac{\|B_{n+1}x\|}{\|B_nx\|}v_{n+1},~(n\geq0).$$ This
shows that ${S_u}_{|{\widehat{\h}}(x)}$ is an injective scalar
unilateral weighted shift with weight sequence
$\big(\frac{\|B_{n+1}x\|}{\|B_nx\|}\big)_{n\geq0}$. Therefore, the
identity, (\ref{ddcc}), follows from theorem 3.7 of
\cite{bourhim}.
\end{proof}

Unlike the scalar weighted shift operators, generally we do not
have $r_1(S_u)=r(S_u)$ if the unilateral operator weighted shift
$S_u$ possesses Bishop's property $(\beta)$ (see example
\ref{exmbishop}). But, of course, if $r_1(S_u)=r(S_u)$, then
either $S_u$ possesses Bishop's property $(\beta)$, or
$\sigma_{\beta}(S_u)=\{\lambda\in\C:|\lambda|=r(S_u)\}$. In
\cite{zguitti}, H. Zguitti represented a unilateral operator
weighted shift as operator multiplication by $z$ on a Hilbert
space of formal power series whose coefficients are in $\h$. He
therefore adapted T. L. Miller and V. G. Miller's arguments given
in \cite{mm} to show that if $S_u$ possesses Bishop's property
$(\beta)$, then $r_2(S_u)=R_1(S_u)$, where
$R_1(S_u)=\liminf\limits_{n\to+\infty}
\big[\inf\limits_{k\geq0}\|B_{n+k}B_k^{-1}\|\big]^{\frac{1}{n}}$.
Here, we refine this result as follows and provide a direct proof.

\begin{thm}
If $S_u$ possesses Bishop's property $(\beta)$, then
$r_2(S_u)=r_1(S_u)$, and $r(S_u)=R_3^+(S_u)$. Moreover, for every
non-zero $x\in\h$, we have
\begin{equation}\label{bbishop}
\lim\limits_{n\to+\infty}
\bigg[\inf\limits_{k\geq0}\frac{\|B_{n+k}x\|}{\|B_kx\|}\bigg]^{\frac{1}{n}}=
\lim\limits_{n\to+\infty}
\bigg[\sup\limits_{k\geq0}\frac{\|B_{n+k}x\|}{\|B_kx\|}\bigg]^{\frac{1}{n}}.
\end{equation}
\end{thm}
\begin{proof}Suppose that $S_u$ possesses Bishop's property $(\beta)$ and
note that, since $S_u$ satisfies Dunford's condition $(C)$,
$r(S_u)=R_3^+(S_u)$ (see theorem \ref{dcc}). If $r_2(S_u)=0$ then,
since $r_1(S_u)\leq r_2(S_u)$, there is nothing to prove. Thus, we
may assume that $0<r_2(S_u)$. Now, recall that it is shown in
\cite{zemanek} that
$$r_1(T)=\min\{|\lambda|:\lambda\in\sigma_{ap}(T)\}$$ for any
operator $T\in\lh$. And so, in order to show that
$r_2(S_u)=r_1(S_u)$, it suffices to prove that ${\mathcal
U}\cap\sigma_{ap}(S_u)=\emptyset$, where ${\mathcal
U}:=\{\lambda\in\C:|\lambda|<r_2(S_u)\}$. Assume for the sake of
contradiction that there is $\lambda_0\in \overline{{\mathcal
U}}\cap\sigma_{ap}(S_u)$. Since $\sigma_p(S_u)=\emptyset$, there
is
$y=(y_0,y_1,y_2,...)\in\cl\big(\ran(S_u-\lambda_0)\big)\backslash\ran
(S_u-\lambda_0)$. For every $x\in\h$, set
$k_x(\lambda):=\sum\limits_{i\geq0}\oplus
\overline{\lambda}^i{B_i^*}^{-1}x,~(\lambda\in{\mathcal U}),$ and
note that
$$(S_u-\lambda)^*k_x(\lambda)=0,~\forall\lambda\in{\mathcal U}.$$
In particular, we have
\begin{equation}
\langle y ,k_x(\lambda_0)\rangle_{\widehat{\h}}=0,~\mbox{ for all
} x\in\h.
\end{equation}
And so, for every $x\in\h$, we have
\begin{eqnarray*}
\langle\sum\limits_{i\geq0}\lambda_0^i{B_i}^{-1}y_i,x\rangle_\h&=&
\sum\limits_{i\geq0}\langle
y_i,\overline{\lambda_0}^i{B_i^*}^{-1}x\rangle_\h\\ &=&\langle y
,k_x(\lambda_0)\rangle_{\widehat{\h}}\\&=&0
\end{eqnarray*}
This implies that
\begin{equation}\label{lambda0}
\sum\limits_{i\geq0}\lambda_0^i{B_i}^{-1}y_i=0.
\end{equation}
Now, for every integer $n\geq0$, we define on ${\mathcal U}$ the
following analytic ${\widehat{\h}}-$valued functions by
$$f(\lambda):=y-\big(\sum\limits_{i\geq0}\lambda^i{B_i}^{-1}y_i\big)^{(0)}
,\mbox{ and
}f_n(\lambda):=y^n-\big(\sum\limits_{i=0}^n\lambda^i{B_i}^{-1}y_i\big)^{(0)},$$
where $y^n:=(y_0,...,y_n,0,0,...)$. Note that for every integer
$n\geq0$, we have $$f_n(\lambda)=\sum\limits_{i=0}^n\big(S_u^i
-\lambda^i\big)\big(B_i^{-1}y_i\big)^{(0)},~(\lambda\in {\mathcal
U}).$$ This implies that each $f_n$ is in
$\ran({(S_u)}_{_{\mathcal U}})$. But
$f\not\in\ran({(S_u)}_{_{\mathcal U}})$ since, in view of
(\ref{lambda0}), we have $f(\lambda_0)=y\not\in\ran
(S_u-\lambda_0)$. On the other hand, for every compact subset $K$
of ${\mathcal U}$, we have
\begin{eqnarray*}
\sup\limits_{\lambda\in
K}\|f_n(\lambda)-f(\lambda)\|_{\widehat{\h}}&\leq&
\|y-y^n\|_{\widehat{\h}}+ \sup\limits_{\lambda\in
K}\|\big(\sum\limits_{i>n}\lambda^i{B_i}^{-1}y_i\big)^{(0)}\|_{\widehat{\h}}\\
&=& \|y-y^n\|_{\widehat{\h}}+ \sup\limits_{\lambda\in
K}\|\sum\limits_{i>n}\lambda^i{B_i}^{-1}y_i\|_\h\\
&\leq&\|y-y^n\|_{\widehat{\h}}+ \sup\limits_{\lambda\in
K}\big\{\sum\limits_{i>n}
|\lambda|^i\|{B_i}^{-1}\|\|y_i\|_\h\big\}\\
&\leq&\bigg(1+\sup\limits_{\lambda\in
K}\big(\sum\limits_{i\geq0}|\lambda|^{2i}\|{B_i}^{-1}\|^2\big)^{\frac{1}{2}}\bigg)\|y-y^n\|_{\widehat{\h}}.
\end{eqnarray*}
Therefore, $f_n\to f$ in ${\mathcal O}({\mathcal
U},{\widehat{\h}})$. As each $f_n\in\ran({(S_u)}_{_{\mathcal U}})$
and $f\not\in\ran({(S_u)}_{_{\mathcal U}})$, we note that
$\ran({(S_u)}_{_{\mathcal U}})$ is not closed. We have a
contradiction to the fact that $S_u$ possesses Bishop's property
$(\beta)$. And so, $\overline{{\mathcal
U}}\cap\sigma_{ap}(S_u)=\emptyset$, as desired.

Now, let $x$ be a non-zero element of $\h$. Since $S_u$ possesses
Bishop's property $(\beta)$, the injective scalar unilateral
weighted shift ${S_u}_{|{\widehat{\h}}(x)}$ possesses also
Bishop's property $(\beta)$. Thus, applying theorem 3.8 of
\cite{bourhim}, gives the identity (\ref{bbishop}).
\end{proof}

\begin{rem} Let $T\in\lh$ be an invertible operator, and assume that $A_n=T$
for all $n\geq0$. The corresponding unilateral operator weighted
shift, $S_u$, satisfies the following identities.
$$r(S_u)=R_3^+(S_u)=r(T)\mbox{ and
}r_1(S_u)=r_2(S_u)=R_2^-(S_u)=r_1(T)=\frac{1}{r(T^{-1})}.$$
Indeed, we clearly have $r(S_u)=r(T)$ and
$r_1(S_u)=r_2(S_u)=r_1(T)=\frac{1}{r(T^{-1})}.$ Since,
$R_3^+(S_u)=\sup\big{\{}r_T(x):x\in\h,~x\not=0\big{\}}$, it
follows from proposition 3.3.14 of \cite{ln} that
$R_3^+(S_u)=r(T)$; therefore, the first identity holds. On the
other hand, we have
\begin{eqnarray*}
R_2^-(S_u)&=&\inf\big{\{}\frac{1}{r_{{T^*}^{-1}}(x)}:x\in\h,~x\not=0\big{\}}\\
&=&\frac{1}{\sup\big{\{}r_{{T^*}^{-1}}(x):x\in\h,~x\not=0\big{\}}}.
\end{eqnarray*}
Again, by proposition 3.3.14 of \cite{ln}, we have
$R_2^-(S_u)=\frac{1}{r(T^{-1})}$; and the second identity follows.
\end{rem}

Assume that $T\in\lh$ is an invertible operator and that $A_n=T$
for all $n\geq0$. So, one may think that the corresponding
unilateral operator weighted shift, $S_u$, satisfies Dunford's
condition $(C)$. It turns out that this is not true in general as
the next example shows.

\begin{exm}\label{kim}Let $(e_n)_{n\in\Z}$ be an orthonormal basis of $\h$,
and let $(\omega_n)_{n\in\Z}$ be a positive two-sided sequence for
which
\begin{itemize}
\item[$(a)$]$0<\inf\limits_{n\in\Z}\omega_n\leq\sup\limits_{n\in\Z}\omega_n<+\infty$.
\item[$(b)$]$\limsup\limits_{n\to+\infty}
\big[\omega_0\omega_1...\omega_{n-1}\big]^{\frac{1}{n}}<\lim\limits_{n\to+\infty}
\big[\sup\limits_{k\geq0}(\omega_k\omega_{k+1}...\omega_{n+k-1})\big]^{\frac{1}{n}}.$
\end{itemize}
Let $T$ be the scalar invertible bilateral weighted shift on $\h$,
defined by $$Te_n=\omega_ne_{n+1},~(n\in\Z).$$ If $A_n=T$ for all
$n\geq0$ then, in view of $(b)$, neither the identity (\ref{ddcc})
nor the identity (\ref{bbishop}) is satisfied for $e_0$. Hence,
$S_u$ is without Dunford's condition $(C)$.

For the construction of a specific example of a positive two-sided
sequence satisfying the above conditions, we refer the reader to
\cite{rid}.
\end{exm}

It is shown in theorem 2.5 of \cite{l.r.williams} that a nonnormal
hyponormal scalar (unilateral or  bilateral) weighted shift has
fat local spectra (see also theorem 3.7 of \cite{bcz}). The next
example shows that this result is not valid for hyponormal
operator weighted shifts.

\begin{exm}\label{exmbishop}Assume that $(e_n)_{n\geq0}$ is an orthonormal basis of $\h$, and
let $(\alpha_n)_{n\geq0}$ be an increasing positive sequence such
that $\lim\limits_{n\to+\infty}\alpha_n=1$. The diagonal operator,
$T$, with the diagonal sequence $(\alpha_n)_{n\geq0}$ (i.e.,
$Te_n=\alpha_ne_n,~\forall n\geq0$) is invertible and satisfies
$r_1(T)=\alpha_0<r(T)=1$. If $A_n=T$ for all $n\geq0$, then the
unilateral operator weighted shift $S_u$ is subnormal. Indeed, let
$\h_n=\h$ for all $n\in\Z$ and let
$$\widetilde{\h}=\sum\limits_{n\in\Z}\oplus\h_n$$ be the Hilbert
space of the two--sided sequences $(x_n)_{n\in\Z}$ such that
$$\|(x_n)_{n\in\Z}\|_{\widetilde{\h}}:=\big(\sum\limits_{n\in\Z}
\|x_n\|_\h^2\big)^{\frac{1}{2}}<+\infty.$$ Let $S_b$ be the
bilateral operator weighted shift defined on $\widetilde{\h}$ by
$$S_b(...,x_{-2},x_{-1},[x_0],x_1,x_2,...)=
(...,Tx_{-2},[Tx_{-1}],Tx_0,Tx_1,...),$$ where for an element
$x=(...,x_{-2},x_{-1},[x_0],x_1,x_2,...)\in\widetilde{\h}$,
$[x_0]$ denotes the central $(0th)$ term of $x$. Note that, since
$T$ is an hermitian operator, $S_b$ is a normal extension of
$S_u$. This shows that $S_u$ is a subnormal operator. Now, we note
that for every $k\geq0$, we have
\begin{eqnarray*}
r_{S_u}(e_k^{(0)})&=&\limsup\limits_{n\to+\infty}\|T^ne_k\|^{\frac{1}{n}}\\
&=&\alpha_k<r(S_u)=1.
\end{eqnarray*}
This shows, on the one hand, that $S_u$ is without fat local
spectra and, on the other hand, that $$r_1(S_u)=r_2(S_u)=
R_2^\mp(S_u)=r_3(S_u)= R_3^-(S_u)=\alpha_0<R_3^+(S_u)= r(S_u)=1.$$
Therefore, in view of the fact that
$\sigma(S_u)=\sigma_{ap}(S_u)\cup\overline{\sigma_p(S_u^*)}$,
corollary \ref{straightforward2} and lemma \ref{spectreponctuel},
we have
$$\sigma_{ap}(S_u)=\{\lambda\in\C:\alpha_0\leq|\lambda|\leq1\}.$$
\end{exm}

\begin{rem} Let $T\in\lh$ be an invertible operator. If $A_n=T$
for all $n\geq0$, then $S_u$ is hyponormal if and only if $T$ is
also hyponormal. Therefore, to construct a kind of example
\ref{exmbishop}, it suffices to take $T$ a hoponormal operator for
which there is a non-zero element $x\in\h$ with $r_T(x)<r(T)$.
\end{rem}

Finally, we would like to point out that

$(a)$ proposition 3.9 of \cite{bourhim} remain valid for the
general setting of operator weighted shift. This is not the case
for proposition 3.11 of \cite{bourhim} as it is shown in example
\ref{kim}.

$(b)$ after the present note was completed, we began to study the
local spectra of bilateral operator weighted shifts; this case is
quite difficult. However, we provided some local spectral
inclusions and obtained a necessary and sufficient condition for a
bilateral operator weighted shift to enjoy the single--valued
extension property. Furthermore, we gave necessary and sufficient
conditions for such operator to satisfy Dunford's condition $(C)$
or Bishop's property $(\beta)$. These results are still on a
preliminary level, and will appear somewhere else once we get some
interesting improvements.

\end{document}